\documentclass[12pt]{amsart}
\usepackage{amsmath, amssymb, amscd, array}

\usepackage{hyperref}

\usepackage{DLdef1}

\newcommand {\del}{{\partial}}

\newcommand {\Lra} {\Longrightarrow}
\newcommand {\Llra} {\Longleftrightarrow}
\newcommand {\wht}{{\mathrm{wht}}}

\let\ssec\subsection

\renewcommand {\ssbegin}[2][*]
 {\refstepcounter{subsection}%
\if#1*
\addcontentsline{toc}{subsection}{\thesubsection.\hskip 1pc #2}%
\else
\addcontentsline{toc}{subsection}{\thesubsection.\hskip 1pc #2. #1}%
\fi
 \def \secno {\gdef \secno {}{\ssecfont
\thesubsection.\hskip 2ex}%
 }%
 \begin{#2}}

\renewcommand {\sssbegin}[2][*]
 {\refstepcounter{subsubsection}
\if#1*
\addcontentsline{toc}{subsubsection}{\thesubsubsection.\hskip 1pc #2}%
\else
\addcontentsline{toc}{subsubsection}{\thesubsubsection.\hskip 1pc #2. #1}
\fi
 \def \secno {\gdef \secno {}{\ssecfont \thesubsubsection.\hskip 2ex}%
 }%
 \begin{#2}}

\renewcommand {\parbegin}[2][*]
 {\refstepcounter{paragraph}
\if#1*
\addcontentsline{toc}{paragraph}{\theparagraph.\hskip 1pc #2}%
\else
\addcontentsline{toc}{paragraph}{\theparagraph.\hskip 1pc #2. #1}
\fi
 \def \secno {\gdef \secno {}{\ssecfont \theparagraph.\hskip 2ex}%
 }%
 \begin{#2}}

\begin{document}

\title{Analogs of Bol operators for $\fpgl(a+1\vert b)\subset \fvect(a\vert b)$}

\author{Sofiane Bouarroudj${}^{a*}$, Dimitry Leites${}^{a,b}$, Irina~Shchepochkina${}^c$}
\address{${}^a$Division of Science and Mathematics\\
New York University Abu Dhabi\\
Po Box 129188, United Arab Emirates\\
${}^{*}$ The corresponding author\\
\email{sofiane.bouarroudj@nyu.edu}\\
${}^b$Department of Mathematics\\
Stockholm University, Stockholm, Sweden\\
\email{dl146@nyu.edu, mleites@math.su.se}\\
${}^c$Independent University of Moscow,
B. Vlasievsky per., d. 11, RU-119 002 Moscow, Russia\\
\email{irina@mccme.ru}\\
}

\begin{abstract}
Bol operators (Bols for short) are differential operators  invariant under the projective action of $\mathfrak{pgl}(2)\simeq\mathfrak{sl}(2)$ between spaces of weighted densities on the 1-dimensional manifold.  

Here, we described analogs of Bols: $\mathfrak{pgl}(a+1\vert b)$-invariant  differential operators between spaces of tensor fields  on $(a\vert b)$-dimensional supermanifolds with irreducible, as $\mathfrak{gl}(a\vert b)$-modules,  fibers of arbitrary, even infinite, dimension for certain ``key" values of $a$ and $b$ --- the ones for which the solution is describable. We discovered many new operators  for $(a|b)=(2|0), (0|3)$ and for the case of $1\vert 1$-dimensional general superstring which looks like a~most natural superization of Bol's result, additional to the cases of  super analogs of Bols  between spaces of weighted densities on the $1\vert n$-dimensional superstrings with a~contact structure we classified in arXiv:2110.10504.

In the case of fibers of dimension $>1$, there are $(a+b-1)$-parameter families of Bols, whereas there are no non-scalar non-zero differential operators between spaces of weighted densities. These two extreme answers justify the selection of cases here. \end{abstract}

\date{}



\maketitle

\markboth{\itshape Sofiane Bouarroudj\textup{,} Dimitry Leites\textup{,} Irina Shchepochkina}{{\itshape Bols for $\fpgl\subset \fvect$}}

\thispagestyle{empty}

\section{Introduction}\label{Sin}

For preliminaries and general overview of Veblen's problem, see \cite{BLS} to which this note is a sequel.
The \textit{Bol operators} (briefly: \textit{Bols}) are unary differential operators   between spaces of weighted densities --- the tensor fields with 1-dimensional fiber --- on the 1-di\-men\-sional domain invariant under the projective action of $\mathfrak{pgl}(2)\simeq\mathfrak{sl}(2)$, see \cite{Bol}.

It is easy to show that on supermanifolds of superdimension $(a|b)\neq (1|0)$, there are no non-scalar non-zero differential operators between spaces of weighted densities.

On the $1\vert n$-dimensional superstrings with a~contact structure, there are  super analogs of Bols  between spaces of weighted densities;  we classified them in \cite{BLS}.

Here, we  consider a maximal  subalgebra $\fpgl(a+1|b)$ in the Lie superalgebra $\fvect(a|b)$ of vector fields with polynomial coefficients on the given superdomain $\cU$ of superdimension $(a|b)$. Here, we classified the following analogs of Bols: $\fpgl(a+1|b)$-invariant operators between spaces of tensor fields  whose fibers are irreducible $\fpgl(a|b)$-modules of any dimension. The answer on the $1|1$-dimensional general superdomains (superstrings) and the cases where $(a|b)=(2|0)$, and $(0|b)$ with infinite-dimensional fibers are the main novelties in the following theorem.

\ssbegin[On $\fpgl(a+1\vert b)$-invariant Bols]{Theorem}[On $\fpgl(a+1|b)$-invariant Bols]\label{Th} \textup{1)} The list of  $\fpgl(a+1|b)$-invariant non-scalar non-zero differential operators between spaces of tensor fields on the superdomain $\cU$ of  ``key'' superdimension $(a|b)$,
is the union of lists in Lemmas $\ref{L02lev1}$, $\ref{L03lev1}$, $\ref{L20}$, $\ref{L11}$. The answers are mostly given in terms of singular vectors, see Subsection~$\ref{ssSing}$.

\textup{2)} For the other values of $(a|b)$, the classification problem of Bols is indescribable, see Subsection~$\ref{wild}$.

\textup{3)} There are no non-scalar non-zero differential $\fpgl(a+1|b)$-invariant operators between spaces of weighted densities on supermanifolds of generic superdimension $(a|b)$, i.e., if ${\dim V=1}$, then $T(V)$ has no $\fpgl(a+1|b)$-submodules.
\end{Theorem}

To double-check our new results we reproduced (for supermanifolds of key dimensions) the classification of differential operators invariant with respect to all changes of variables, cf. \cite{BL, LKW, GLS}. These results are called Claims. In particular, we corrected one claim in  \cite{LKW}, see Claim~\ref{ind}.

The computations in Lemmas and Claims (Theorem~\ref{Th}) are straightforward (although involved in the cases of Lemma~\ref{L20}) and all alike (``reading them will not make the reader wiser''). We collected them, except the proof of Claim~\ref{ind} with its interesting answer, in \S~\ref{Sbor} and skipped the easy proof of item 3) of  Theorem~\ref{Th}.

\section{Bols invariant under $\fpgl(a+1|b)\simeq\fpgl(b|a+1)\subset \fvect(a|b)$}
\label{Sbinv}

\ssec{Restricted and unrestricted duals}
Let $\fg:=\fvect(a|b)$ be the Lie superalgebra of vector fields with polynomial coefficients in indeterminates $x=(u,\xi)$, where $u=(u_1,\dots, u_{a})$ are even indeterminates, and $\xi=(\xi_1, \dots, \xi_{b})$   are odd. Consider
the standard $\Zee$-grading of $\fg$ given by $\deg x_i=1$ for all $i$. Then, $\fg_0\simeq\fgl(a|b)$. Set
\[
\fg_\geq:=\oplus_{i\geq 0}\fg_i, \ \ \ \fg_>:=\oplus_{i> 0}\fg_i; \ \ \Zee_{> 0}:=\{1, 2, \dots \}.
\]

Let $V$ be an irreducible $\fg_0$-module with either highest or lowest weight vector (or with both, if $V$ is finite-dimensional). Make $V$ a~$\fg_\geq$-module by setting $\fg_>V=0$.Observe that $\fg =\fg_{-1} \oplus \fg_\geq$ and the Lie superalgebra $\fg_{-1}$ is commutative; set (P.B.W. is for the Poincar\'e--Birkhoff--Witt theorem)
\[
I(V):=U(\fg)\otimes_{U(\fg_\geq)}V\simeq U(\fg_{-1})\otimes U(\fg_\geq)\otimes_{U(\fg_\geq)}V\stackrel{\text{P.B.W.}}{\simeq}\Cee[\fg_{-1}]\otimes V=\Cee[\partial_x]\otimes V.
\]

Recall the following isomorphisms of the dual spaces $(\Cee[x])^*\simeq \Cee[[x^*]]$ and $(\Cee[[x]])^*\simeq \Cee[x^*]$, where $x^*$ is the vector dual to $x$.  We will be working with induced $\fg$-modules $I(V)$ (whose elements are polynomials) but interpret our results in terms of spaces $T(V^*)$ of tensor fields which are $V^*$-valued formal power series dual to the elements of $I(V)$:
\[
T(V^*):=\Hom_{U(\fg_\geq)}(U(\fg), V^*)\stackrel{\text{P.B.W.}}{\simeq}\Cee[[(\fg_{-1})^*]]\otimes V^*=\Cee[[x]]\otimes V^*.
\]

We will also need the dual spaces of tensor fields expressed in the same terms (for example, spaces of integrable forms dual to the spaces of differential forms). Researchers who, for whatever reason, would like to express the dual of the space of polynomials as a~ space of polynomials, not power series, introduced \textit{restricted dual} $W^\#:=\oplus_{i\in \Zee}(W_i)^*$ of the graded space $W=\oplus_{i\in \Zee}W_i$. For example, $(\Cee[x])^\#\simeq \Cee[x^*]$, where we can identify $x^*=\partial_x$ and $(\partial_x)^*=x$. 

Likewise, working with duals of invariant operators we would like, for the sake of uniformity, to express the duals of the space of series as a~ space of series, not polynomials. We define the \textit{unrestricted dual} of $\Cee[[x]]$ by setting
\[
(\Cee[[x]])^\flat:= \Cee[[x]]\vvol,
\]
where $\vvol:=\vvol(x)$ is the volume element in coordinates $x$ such that the change of coordinates $x\mapsto y(x)$ transforms volume elements as follows (on supermanifolds, we take for $\det$ the superdeterminant a.k.a. Berezinian)
\[
f(y)\vvol(y) \mapsto f(y(x))\cdot \det\left(\nfrac{\partial y}{\partial x}\right)\cdot\vvol(y(x)).
\]
The origin of what we called \textit{unrestricted dual} --- we did not encounter a~ special term for it in the literature --- is as follows.

Consider the case of supermanifolds of dimension $a|b$. Let $\cF$ denote the algebra of functions of the type we are considering (smooth, polynomial, formal power series, etc.).
Observe that $T(V)\otimes_\cF T(W)\simeq T(V\otimes W)$. Consider the paring given by integration of smooth  functions with compact support, or rapidly decaying at infinity, i.e., consider \textbf{dualization over the algebra~ $\cF$}. For such functions, we see that
\be\label{approx}
(T(X))^\flat\approx T(X^*)\otimes_\cF \Vol\simeq T(X^*\otimes \str) ,
\ee
where  $X$ is any (irreducible) $\fgl(a|b)$-module,
$\Vol$ is the space of densities a.k.a. volume forms, $\str$ denotes both the supertrace and a~1-dimensional $\fgl(a|b)$-module given by the supertrace, so $\Vol=T(\str)$. For example, if $X=\Span(\partial_{x_1}, \dots ,\partial_{x_m})$ for $m=a+b$, then $X^*=\Span(dx_1, \dots , dx_m)$, each differential $dx_i$ being considered as a~ \textbf{linear functional over~ $\cF$ on the $\cF$-module of vector fields}:
\[
fdx_i \ (g\partial_{x_j}):=(-1)^{p(g)(p(x_i)+1)}\delta_{ij}\ fg,\text{~~where $f,g\in\cF$}.
\]

If elements of one of the  spaces $(T(X))^*$ or $T(X)$ are smooth (or generalized) vector-valued functions with compact support (or rapidly decaying at infinity or whatever allowing us to integrate the ``convolution'' $t(s)\in\cF$, where $t\in T(X^*)$ and $s\in T(X)$), then, recalling that one can only integrate volume forms, not functions, we see that  $T(X^*)\otimes_\cF \Vol$ and $T(X)$ are manifestly in duality given by the integral. 

For example, on the $n$-dimensional manifold, we would like to consider the dual of $d:\cF=\Omega^0\tto\Omega^1$, where $\cF$ is the space of functions,  as the operator $d:\Omega^{n-1}\tto\Omega^n=\Vol$ or, more generally, to consider the dual of $d_i:\Omega^i\tto\Omega^{i+1}$ as $d_{n-1-i}:\Omega^{n-i-1}\tto\Omega^{n-i}$. All these differentials are denoted by the same symbol $d$ (since all of them are of the form $\sum \widehat x_i\partial_{x_i}$, where $ \widehat x_i:=dx_i$) and called the \textit{exterior differential}.

In our case of polynomial or formal coefficients, there is no integration, we just interpret the above isomorphism $\approx$  in \eqref{approx} in a~ ``Pickwickian sense".
Namely, in precise terms, the sign $\approx$ in \eqref{approx} is an isomorphism ($\simeq$) if the spaces on both sides of $\approx$ are of finite dimension; otherwise it is the definition ($:=$) of the unrestricted dual.

\ssec{On invariant operators and singular vectors}\label{ssSing} To every invariant differential operator $D:T(X)\tto T(Y)$, where $X$ and $Y$ are $\fg_0$-modules with lowest weight vectors, there corresponds the dual over the ground field operator $D^*: I(Y^*)\tto I(X^*)$ since $(T(X))^*\simeq I(X^*)$ and $(I(X^*))^*\simeq T(X)$. Then, to a~non-zero $D$, there corresponds a~submodule in $I(X^*)$; the $\fg_0$-highest weight vector of this submodule is called \textit{singular}. So, a~non-zero vector $f_i\in I(X^*)_{-i}$ of degree $i$ in a~ basis of $\fg_{-1}$ --- we will also say that $i$ is the \textit{level} of $f_i$ --- is singular if and only if it is $\fg_0$-highest weight vector, and $\fg_>(f_i)=0$.

\ssec{Remarks} 1) We use Rudakov's approach, see \cite{R1} and elucidations in \cite{GLS}, and seek singular vectors in the induced modules $I(V)$ because it is much easier than to directly describe matrix-valued differential operators between spaces of tensor fields $T(V)$, although we will interpret our singular vectors in terms of operators between spaces $T(V)$.  The direct approach (computing in terms of $T(V)$) is feasible only if $\dim V=1$ or 2.

2) \textbf{Over any manifold (domain, since the problem is local, see \cite{BL}) $\cU$, the singular vectors in $I(V)$ exist only if $\dim V<\infty$; there are only $\dim\, \cU$ possible irreducible such modules, namely $V=E^i(\id)$ corresponding to $d:\Omega^i\tto \Omega^{i+1}$}, see \cite{R1} and its clarification in \cite{BL} and \cite{GLS}. On superdomains of any dimension ${a|b}$ with $b\neq 0$, there are $\fpgl(a+1|b)$-singular vectors in $I(V)$ even if $\dim V=\infty$ and there are ``too many'' of them, except for the ``key dimensions" $(a|b)$ considered in Lemmas below.

3) As abstract Lie superalgebras, $\Zee$-grading ignored,
\be\label{iso}
\text{$\fvect(1|1)\simeq \fk(1|2)$ and $\fsl(2|1)\simeq \fosp(2|2)$.}
\ee
Hence, it looks like ``since we have described Bols in terms of $\fk(1|2)$, see \cite{BLS}, we have already got  Bols in the $\fsl(1|2)\subset \fvect(1|1)$ case". But this is not so!

 Indeed, the isomorphisms \eqref{iso} consider algebras as abstract, whereas we are interested in geometry preserved by the filtered algebras and their graded subalgebras. In particular, the spaces of weighted densities over $\fvect(a|b)$ depend on one parameter, whereas the spaces of weighted densities over $\fk(1|2)$ depend on two parameters.
For $\fsl(1|2)\subset \fvect(1|1)$, we describe Bols  \textit{ab ovo}, see Lemma~\ref{L11}.

4) Let $\fg\neq\fvect(1|0)$, $\fk(1|1)$ or $\fk(1|2)$. Then, $\fg_>$ is generated by $\fg_1$ and it suffices to consider only the conditions $Lv=0$ for the $\fg_0$-lowest  weight operators $L\in\fg_1$,  instead of $\fg_>v=0$.

In the cases of $\fg\simeq\fvect(1|0)$, $\fk(1|1)$ or $\fk(1|2)$, in order to describe $\fg$-singular vectors, we have to consider the lowest weight vectors of $\fg_2$. This is not needed to describe Bols.

5) \textbf{If $V$ has no lowest weight vector, then $V^*$ has no highest weight vector and $I(V^*)$ can have no singular vectors, although the dual operator exists}. This fact was first formulated in \cite{GLS}. Thus, there are invariant operators $D: T(V)\tto T(W)$ that can not be discovered by searching for singular vectors.

\ssec{Differential and integrable forms}\label{DifInt} Let $\cF:=\Cee[u,\xi]$ be the space of functions on an $a|b$-dimensional supermanifold. Then, the space
$\Omega^\ell$ of differential $\ell$-forms is the $\cF$-module with the basis of homogeneous degree-$\ell$ polynomials in $\widehat x_i:=dx_i$, where $p(\widehat x_i)=p(x_i)+\od$  for  $x=(u, \xi)$, while the space
$\Sigma_{-k}$ of integrable $(-k)$-forms is the $\cF$-module with the basis of homogeneous degree-$k$ polynomials in $\del_{\widehat x}$  with values in the space  $\Vol(x):=\cF\vvol(x)$ of volume forms, where $\vvol$ is the volume element with constant coefficients, i.e., 
\[
\text{the class of the element ${du_1\wedge\dots \wedge du_a\partial_{\xi_1}\dots \partial_{\xi_b}}$}
\]
 in the quotient of (recall that $E^k(V):=S^k(\Pi(V))$ for any superspace $V$, where $\Pi$ is the change of parity functor, $E^k$ is the operator of raising to the $k$th \textit{exterior} power and $S^k$ is the operator of raising to the $k$th \textit{symmetric} power; the index indicates the algebra over which these operators are considered if other than the ground field)
\be\label{pseu}
E^a_\cF((\fvect(a|b))^*)\otimes_\cF S^b_\cF(\fvect(a|b))
\ee
--- the tensor product over $\cF$ of the $a$th exterior power of the $\cF$-module of covector fields by the $b$th symmetric power of the $\cF$-module of vector fields --- modulo the maximal $\fvect(a|b)$-submodule of corank 1. In \cite{BLi}, the elements of $\Sigma_{\bcdot}$ were called \textit{integrable} forms because they can be integrated; Deligne calls them \textit{integral}, see \cite[Ch.~3]{Del}, either as in the (incorrect, we think) translation from the Russian (\cite{BLi}) or to rhyme with differential.

\sssec{Pseudo-differential and pseudo-integrable forms on supermanifolds of dimension $(a|1)$}\label{PseuA1} Recall (see \cite{BL,GLS}) that the case where $\sdim\, \cU=(a|1)$ is special from the point of view of $\fvect(a|1)$-modules $T(V)$.
In this case, we can let the degree of the differential or integrable forms be any $\lambda\in\Cee$ and the dimension of the fiber of the corresponding tensor fields remains finite.

Indeed, set 
\[
\text{$\widehat \xi^{-1}:= \del_{\widehat \xi}$, $\widehat \xi^0:=1$  and $\widehat \xi\cdot \widehat \xi^\lambda:=\widehat \xi^{\lambda+1}$ for any $\lambda\in\Cee$; let $\deg(\widehat \xi^\lambda):=\lambda$.}
\]
 Set $\Psi\Omega^{\bcdot}:=\oplus_{\lambda\in\Cee}\Psi\Omega^\lambda$, where $\Psi\Omega^\lambda$ is the $\Cee[[u,\xi]]$-module with basis
\[ \widehat \xi^\lambda,\  \widehat \xi^{\lambda-1}\widehat u_i, \ \widehat \xi^{\lambda-2}\widehat u_i\widehat u_j, \dots ,
\widehat \xi^{\lambda-a}\widehat u_1\widehat u_2\dots \widehat u_a.
\]
Let $\Psi\Omega_\Zee^{\bcdot}:=\oplus_{\lambda\in\Zee}\Psi\Omega^\lambda$. Define the maps $\alpha:\Omega^{\bcdot}\tto \Psi\Omega_\Zee^{\bcdot}$ and $\beta:\Psi\Omega_\Zee^{\bcdot}\tto\Sigma_{\bcdot}$ by setting
\[
\alpha(\omega)=\omega \text{~~and ~~}\beta(f(u,\xi)\widehat \xi^{k}\widehat u_1\widehat u_2\dots \widehat u_a)=f(u,\xi)\del_{\widehat \xi}{}^{-k+1}\vvol(u,\xi),
\]
since $\vvol:=\vvol(u,\xi)$ is the class of $\del_{\xi}\widehat u_1\widehat u_2\dots \widehat u_a$ on which $\fvect(a|1)$ acts as on $\del_{\widehat \xi}\widehat u_1\widehat u_2\dots \widehat u_a$. Clearly, the action of the exterior differential can be extended from the space of forms to the space $\Psi\Omega^{\bcdot}$ of \textit{pseudo-differential} forms of complex degrees. It can also be extended from the space $\Sigma_{\bcdot}$ of  integrable forms to the space $\Psi\Sigma_{\bcdot}$ of \textit{pseudo-integrable} forms, where $\Psi\Sigma_{\lambda}$ is the the $\Cee[[u,\xi]]$-module with basis
\[ \widehat \xi^\lambda\vvol,\ \ \widehat \xi^{\lambda-1}\del_{\widehat u_i}\vvol, \ \ \widehat \xi^{\lambda-2}\del_{\widehat u_i}\del_{\widehat u_j}\vvol, \dots ,
\widehat \xi^{\lambda-a}\del_{\widehat u_1}\del_{\widehat u_2}\dots \del_{\widehat u_a}\vvol.
\]

\ssbegin{Lemma}[\cite{BL, LKW}]\label{L02lev1} Let $\sdim\, \cU=(0|2)$. In this case, there are as many Bols as there are $\fvect(0|2)$-invariant differential operators because $\fsl(1|2)\simeq\fvect(0|2)$.

\textup{1)} If $\dim V<\infty$, then the degree-$1$ invariant operator is the exterior differential of one of the types $d:\Omega^\ell\tto \Omega^{\ell+1}$ or
$d:\Sigma_{-k}\tto\Sigma_{-k+1}$, where $k\in\Zee_{>0}$ and $\ell\in\Zee_{\geq 0}$.

\textup{2)}  If $\dim V=\infty$, then  $k\in\Cee\setminus\Zee_{>0}$ and $\ell\in\Cee\setminus\Zee_{\geq 0}$.  The corresponding spaces $I(V)$ are dual to the spaces of pseudo-differential forms of degree $\ell$ and pseudo-integrable forms of degree $k$. In both cases, the invariant operator is the exterior differential.

The $\fg$-invariant degree-$2$ operator is the Berezin integral. \end{Lemma}

\ssec{$\sdim \cU=(0|b)$} For $b>2$, there should be more Bols than
there are $\fvect(0|b)$-invariant operators, see \cite{LKW} and Subsection~\ref{ind}.

\sssbegin{Lemma}[The  key case where $b=3$]\label{L03lev1} Let $\sdim\, \cU=(0|3)$. Only the $\fsl(1|3)$-invariant Bols corresponding to the following singular vectors are possible. 

\begin{tabular}{|c|c|c|}
\hline
\multicolumn{3}{|c|}{\bf level 1: $f_1=\del_1v_1+\del_2v_2+\del_3v_3$}\\
\hline
Case & Highest weight $\wht$ of $V$& 
Conditions on $\wht$ if $\dim V<\infty$\\
\hline
$v_1=v_2=0$, $v_3\ne 0$ & $(\lambda, -\lambda, \nu)$ &$\nu=-\lambda-k$, where $\lambda, k\in \Zee_{\geq0}$\\
\hline
$v_1=0$, $v_2\ne 0$ & $(\lambda, \mu,-\lambda+1)$ &$\wht(v_2)=\left(\nfrac n2, \nfrac n2-k, -\nfrac n2+1\right)$,\\
&& where  $n\in\Zee_{> 0},\ \  k=0,1,2,\dots, n-1$\\
\hline
$v_1\ne 0$ & $(\lambda, \mu,-\mu+2)$ &$\wht(v_1)=\left (\nfrac n2+k, \nfrac n2+1, -\nfrac n2+1\right)$,\\
&& where $n\in\Zee_{\ge 0},\ \  k\in\Zee_{> 0}$  \\
\hline
\multicolumn{3}{|c|}{\bf level 2: $f_2=\del_2\del_3v_1+\del_2\del_1v_2+\del_1\del_2v_3$}\\
\hline
$v_2=v_3=0$, $v_1\ne 0$ & $(\lambda, -\lambda+1, -\lambda+1)$ & $\wht(v_1)=\left(\nfrac n2,-\nfrac n2+1, -\nfrac n2+1\right)$,
\\
&& where $n\in\Zee_{> 0}$\\
\hline
$v_3=0$, $v_2\ne 0$ & $(\lambda, -\lambda+1,\lambda+1)$ &no solutions\\
\hline
$v_3\ne 0$ & $(\lambda, \lambda,-\lambda+2)$ &$\wht(v_3)=\left (\nfrac n2+1,\nfrac n2+1,-\nfrac n2+1\right)$, \\
&& where $n\in\Zee_{\ge 0}$\\
\hline
\multicolumn{3}{|c|}{\bf level 3: $f_3=\del_1\del_2\del_3v$. Then, $\wht(v)=(1,1,1)$ and $\wht(f_3)=(0,0,0)$}\\
\hline
\end{tabular}
\end{Lemma}

\sssbegin{Claim}[\cite{BL,LKW}] \label{2.3.2} \label{Cl1}
Let $f_1=\del_1v_1+\del_2v_2+\del_3v_3\in I(V)$ be a $\fvect(0|3)$-singular vector. Then, the following cases can occur.

\textup{1)} If $v_1=v_2=0$, we have
\[
\wht (v_3)=(0,0,-k) \text{~~and $\wht (f_1)=(0,0,-k-1)$, where $k\in\begin{cases}\Zee_{\geq0}&\text{if $\dim V<\infty$},\\
\Cee\setminus\Zee_{\ge 0}&\text{if $\dim V=\infty$}.\end{cases}$}
\]
This corresponds to the exterior differential between spaces of (pseudo)integrable forms.

\textup{2)} If $v_1=0, v_2\ne 0$, we have
\[
\begin{cases}\text{no solutions}&\text{if $\dim V<\infty$}\\
\wht(v_2)=(0,\mu,1)\text{~~and $\wht (f_1)=(0,\mu-1,1)$},&\\
\text{where $\mu\in\Cee$},&\text{if $\dim V=\infty$}.\end{cases}
\]
This corresponds to the exterior differential between spaces of pseudo-forms with $\fgl(3)$-lowest weight vectors of the form $\widehat \xi_1(\widehat \xi_2)^{\mu}\vvol(\xi)$.

\textup{3)} If $v_1\ne 0$, we have
\[
\wht(v_1)=(m,1,1) \text{~~and $\wht (f_1)=(m-1,1, 1)$, where $m\in\begin{cases}\Zee_{>0}&\text{if $\dim V<\infty$},\\
\Cee\setminus\Zee_{> 0}&\text{if $\dim V=\infty$}.\end{cases}$}
\]
This corresponds to the exterior differential between spaces of (pseudo)differential forms.

There are no singular vectors of level $2$. 

The singular vectors of level $3$ corresponds to the Berezin integral.
\end{Claim}

\sssbegin{Claim}[\cite{LKW}]\label{ind}
Let $f_1=\sum\del_iv_i\in I(V)$ be a~$\fvect(0|b)$-singular vector. Then, it corresponds to the exterior differential $d$ and the highest weight of $V$ can be only of the form
$(0, \dots , 0, \underbrace{\lambda}_r, 1, \dots , 1)$,  where
\[
\begin{array}{ll}
\text{if~~}r=1, &\text{~~then}\begin{cases}
\lambda\in\Zee_{>0}\text{ and }&\text{if $\dim V<\infty$},\\
\lambda\in \Cee\setminus \Zee_{> 0}&\text{if $\dim V=\infty$};\end{cases}\\
 &\text{corresponds to $d$ between spaces of (pseudo)integrable forms;}\\
\text{if~~}r=b, &\text{~~then}\begin{cases}\lambda\in\Zee_{\leq 0}&\text{if $\dim V<\infty$},\\
\lambda\in \Cee\setminus \Zee_{\leq 0}&\text{if $\dim V=\infty$};
\end{cases}\\
&\text{corresponds to $d$  between spaces of (pseudo)differential forms;}\\
\text{if~~}1<r<b, &\text{~~then $\lambda\in \Cee$ and  $\dim V=\infty$}\\
 &\text{corresponds to $d$   between spaces of pseudo-forms}\\
 &\text{with $\fgl(b)$-lowest weight vectors of the form $\widehat \xi_1\dots \widehat \xi_{r-1}(\widehat \xi_r)^{\lambda+1}\vvol(\xi)$}.
\end{array}
\]

For $f_i$ where $1<i<b$, there are no singular vectors \textup{(contrary to \cite[item 1.2, p.225]{LKW})}.

For $f_b$, the highest weight of $V$ is $(1, \dots , 1)$ and
$\wht(f_b)=(0,\dots, 0)$; corresponds to the Berezin integral.
\end{Claim}

\textbf{Comment}. Let us explain: Since $\vvol=\prod \partial_{\xi_i}$, then, from the point of view of the $\fvect(0|b)$-action, the product $\widehat \xi_1\dots \widehat \xi_{b-1}(\widehat \xi_b)^{s}\vvol(\xi)$ is  the same as $(\widehat \xi_b)^{s-1}$ which is the lowest weight vector in the $\fgl(b)$-module of (pseudo)differential forms with constant coefficients. If $r=1$, recall that $(\widehat \xi_1)^{s}=(\partial_{\widehat \xi_1})^{-s}$. Since this claim is published only in a~ not easily accessible collection \cite{LKW}, and the answer is not obvious, we give its proof for the sake of completeness and because the induction base in \cite{LKW} contains a mistake (due to Leites); for a correction, see Claim~ \ref{2.3.2}.

\begin{proof} We proceed by induction. \textbf{The induction base} is proved for $b=3$, see Claim~ \ref{2.3.2}.

\textbf{The induction step}: In $\fvect(0|b+1)=\fvect(\xi_1,\dots,\xi_{b+1})$, we consider  two subalgebras $\fg=\fvect(\xi_1,\dots,\xi_b)$ and $\fh=\fvect(\xi_2,\dots,\xi_{b+1})$.

Set
\[
\begin{array}{l}
\fg^+:=\Span(\xi_i\del_j\mid i<j\le b)\oplus \fg_{>0},\\
\fh^+:=\Span(\xi_i\del_j\mid 2\leq i<j\le b+1)\oplus \fh_{>0}.
\end{array}
\]
 Clearly,
\[
[\fg^+,\del_{b+1}]=0, \quad [\fh^+,\del_{1}]=0.
\]

For any $k<b+1$, we have (this is a definition of the vectors $g_i$ and $h_i$)
\[
f_k=g_k+\del_{b+1}g_{k-1}=h_k+\del_1 h_{k-1}.
\]

Any $\fvect(0|b+1)$-singular vector is simultaneously $\fg$-singular and $\fh$-singular. Therefore, if $f_k$ is singular, then
\[
Y(f_k)=Y(g_k)+\del_{b+1}Y(g_{k-1})=0 \text{ for all } Y\in \fg^+,
\]
implying $Y(g_k)=0$, i.e., $g_k$ is $\fg$-singular. Hence, if $1<k<b+1$, then $g_k=0$. Similarly, $h_k=0$. This means that there are no singular vectors on level $k$, where $1<k<b+1$.

Now, let $k=1$. 
Let $f_1=\del_1 v_1+\del_2 v_2+\dots +\del_b v_b$.

Set $X_{ij}^+:=\xi_i\del_j$, where $i<j$. Then, $[X_{ij}^+,\del_r]=-\delta_{ir}\del_j$, and hence
\[
X_{ij}^+(f_1)=\sum_{1\leq r\leq b} \del_r(X_{ij}^+ v_r)-\del_j v_i.
\]
Therefore,
\[
X_{ij}^+(f_1)=0 \Llra \begin{cases}
                            X_{ij}^+v_r=0 &\text{for $r\ne j$},\\
                            X_{ij}^+v_j=v_i.\\
                            \end{cases}
 \]

We see that if $v_1=\dots=v_{r-1}=0$, and $v_r\ne 0$, then $v_r$ is the highest weight vector. Let us see how the elements of $\fg_1$ act on it:
\[
[\xi_r\xi_i\del_i,\del_r]=H_i,\; [\xi_r\xi_i\del_i,\del_i]=-(\xi_r\del_i)\Lra \xi_r\xi_i\del_i(f_1)=H_i v_r-(\xi_r\del_i) v_i.
\]

This immediately implies that for $i<r$, the condition $\xi_r\xi_i\del_i(f_1)=0$ is equivalent to  $H_i v_r=0$, whereas if $i>r$, then $\xi_r\del_i=X_{ri}^+$ and $X_{ri}^+v_i=v_r$, and hence the condition $\xi_r\xi_i\del_i(f)=0$ for $i>r$ is equivalent to $H_iv_r=v_r$. Thus, for any $\lambda\in\Cee$, we have
\[
\wht(v_r)=(0,\dots, 0, \underbrace{\lambda}_{r}, 1, \dots, 1),\quad \wht(f_1)=(0,\dots, 0, \underbrace{\lambda-1}_{r}, 1, \dots, 1).
\]

Finally, $f_{b+1}=\del_{b+1}g_b=\del_1h_b$, and hence if $f_{b+1}$ is singular, then $g_b$ and $h_b$ are singular. By the induction hypothesis
\[
\wht(g_b)=(0,\dots,0,\nu)\Lra \wht(f_{b+1})=(0,\dots,0,\nu-1) \text{ and }
\]
\[
\wht(h_b)=(\nu',0,\dots,0)\Lra \wht(f_{b+1})=(\nu'-1,0,\dots,0).
\]
Comparing these expressions we see that $\nu=\nu'=1$ and $\wht(f_{b+1})=(0,\dots,0)$. \end{proof}

\ssec{The case of $\sdim \cU=(a|0)$} The $\fvect(a|0)$-invariant operators are described by Rudakov \cite{R1}, see review \cite{GLS}, where Rudakov's results are interpreted: no non-zero non-scalar differential operator except for the exterior differential $d: \Omega^i\tto \Omega^{i+1}$ in the de~Rham complex.

\sssec{The exceptional case of $\sdim \cU=(1|0)$} This case was considered by Bol,  see \cite{Bol}; for an interpretation of Bols in terms of Verma modules over $\fsl(2)$, see \cite{BLS}. 

\sssbegin{Lemma}[The key case of $\sdim \cU=(2|0)$]\label{L20} The following are the only possible singular $\fsl(3)$-invariant vectors $f_n=\sum_{m=0}^n\del_1^m\del_2^{n-m}v_m$, where $X^+v_n =0$, and $X^+v_m = (m+1)v_{m+1}$ for $m=  0,1,\dots, n-1$. There are $n+1$ cases:
\[
\textbf{Case $k$: } v_n=\dots = v_{k+1}=0, \ \ v_k\ne 0, \text{ where } k= 0,1,\dots, n.
\]
In these cases, we have
\[
\textbf{Case $0$: } \wht (v_0)=(n-1-2\lambda, \lambda), \ \ \wht (f_n)=(n-1-2\lambda, \lambda-n),\]
where
\[
\begin{cases}n-1-3\lambda\in\Zee_{\geq 0}&\text{if $\dim V<\infty$},\\
n-1-3\lambda\in\Cee\setminus\Zee_{\geq 0}&\text{otherwise}.
\end{cases}
\]

{\bf Case $k$, where $0< k<n$ and $n\ge 2$}. 
Then, $\dim V=\infty$ and
\[
\wht(v_k)=(\nfrac{2k-n}3-1,  \nfrac{2n-k}3), \ \ \wht(f_n)=(-\nfrac{k+n}3-1, \nfrac{2k-n}3).
\]

{\bf Case $n$}. Then, 
\[
\wht(v_n)=(\lambda, n-2-2\lambda), \ \ \wht(f_n)=(\lambda-n, n-2-2\lambda),
\]
where
\[
\begin{cases}3\lambda-n+2\in\Zee_{\geq 0}&\text{if $\dim V<\infty$},\\
3\lambda-n+2\in\Cee\setminus\Zee_{\geq 0}&\text{otherwise}.
\end{cases}
\]

\end{Lemma}

\ssec{The key case of $\sdim \cU=(1|1)$}{}~{}

\sssbegin{Lemma}\label{L11} Let $\sdim\, \cU=(1|1)$ with coordinates $x$ (even) and $\xi$ (odd). For a~basis of $\fg_{-1}$ we take $\del:=\del_x$ and $\delta:=\del_\xi$.

There are two degree-$n$ Bols for any $\lambda\in\Cee$ corresponding to singular $\fsl(2|1)$-invariant vectors $f_n=\del^n v+\del^{n-1}\delta w\in I(V)$:
\[
\begin{array}{ll}
1) \ v=0,&\text{then~} \wht(w)=(n-1, \lambda),\ \ \wht(f_n)=(0, \lambda-1);\\
2) \ v\neq 0,&\text{then~}  \wht(v)=(\lambda, n-2\lambda),\ \ \wht(f_n)=(\lambda-n, n-2\lambda).\\
\end{array}
\]
\end{Lemma}

\sssbegin{Claim}[\cite{BL}]\label{C11} Let $\sdim\, \cU=(1|1)$. For a~basis of $\fg_{-1}$ we take $\del:=\del_x$ and $\delta:=\del_\xi$.

\textup{1)} The only non-zero non-scalar $\fvect(1|1)$-invariant differential operators of degree $1$ correspond to the following singular vectors $f_1=\del v+\delta w\in I(V)$, where $\lambda\in\Cee$:
\[
\begin{array}{ll}
1) \ v=0,&\text{then~} \wht(w)=(0, \lambda),\ \ \wht(f_1)=(0, \lambda-1);\\
2) \ v\neq 0,&\text{then~}  \wht(v)=(0, 1),\ \ \wht(f_1)=(-1, 1).\\
\end{array}
\]

\textup{2)} There are no non-zero non-scalar $\fvect(1|1)$-invariant differential operators of degree~ $>1$.
\end{Claim}

\sssec{Indescribably many Bols for $\fpgl(a+1|b)\subset \fvect(a|b)$ for $(a|b)$ generic}\label{wild}{}~{} Let $\fvect(a|b)=\fvect(x_1, \dots, x_m)$, where $m=a+b$,  $p(x_i)=\ev$ for $1\le i\le a$  and $p(x_i)=\od$ for  $a+1\le i\le m$. Let
\[
\begin{array}{l}
X_i^+=x_i\del_{i+1} \;  \text{ for } i=1,2,\dots, m-1;\ H_i=x_i\del_i  \text{ for }  i=1,\dots, m;
\\
s_m=x_m\sum_{1\leq i \leq m-1}x_i\del_i. 
\end{array}
\]

Let us consider just one case where $f=\del_m v_m$: already in this case, there manifestly indescribably many $\fpgl(a+1|b)$-singular vectors.
We have
\[
\begin{array}{l}
{}[X_i^+,\del_m]=0 \text{ for all } i=1,\dots, b-1;
\\
{}[s_m,\del_m]=\sum_{1\leq i \leq m-1}H_i;
\\
\end{array}
\]

Then, $X_i^+v_m=0$ for all $i=1, \dots, m-1$, i.e., $v_m$ is the highest weight vector and the condition $s_mf=0$ implies $(H_1+\dots + H_{m-1})v_m=0$.

If $\wht(v_m)=(\lambda_1, \dots, \lambda_m)$, then being a~ Bol imposes only one condition on $\wht(v_m)$:
\be\label{eqwild}
\sum_{1\leq i \leq m-1}\lambda_i=0.
\ee
The additional conditions are imposed by dimension of $V$:
\[
\begin{cases}\lambda_i-\lambda_{i+1}\not\in\Zee_{\ge 0}\text{~for at least one $i$}&\text{if $\dim V=\infty$},\\
\lambda_i-\lambda_{i+1}\in\Zee_{\ge 0}\text{~for all $i$}&\text{if $\dim V<\infty$}.
\end{cases}
\]

\textbf{Thus, even in this particular case, we have $a+b-1$ parameters. And there are $a+b-2$ other cases on level $1$ as well as singular vectors on higher levels. To describe them all does not seem feasible.}

\section{Straightforward calculations}\label{Sbor}

We will use a well-known fact: if the irreducible $\fgl(n)$-module $V$ with highest weight $(\lambda_1, \dots, \lambda_n)$ with respect to the elements $E_{ii}\in\fgl(n)$ is finite-dimensional, then $\lambda_i-\lambda_{i+1}\in\Zee_{\geq0}$ for all $i$ .

\ssec{Proof of Lemma $\ref{L02lev1}$}  Consider the following shorthand of basis elements of $\fvect(0|2)$:
\[
\begin{array}{l}
X^+:=\xi_1\del_2,\ \ X^-:=\xi_2\del_1, \ \ H_i:=\xi_i\del_i; \\
s_1:=\xi_1\xi_2\del_2,\ \ s_2:=\xi_1\xi_2\del_1.
\end{array}
\]

\underline{\textbf{Level 1}: $f_1=\del_1v_1+\del_2v_2\in I(V)$.}
The conditions
\[
X^+f_1=0,\ \ s_1f_1=0,\ \ s_2f_1=0
\]
are equivalent, respectively, to
\[
\begin{array}{l}
X^+v_2=v_1, \ \ X^+v_1=0; \\
X^-v_1-H_1v_2=X^-X^+v_2-H_1v_2=0;\\
H_2v_1-X^+v_2=0.
\end{array}
\]
We arrive at the two cases:
\be\label{e02}
\begin{array}{ll}
1) \ v_1=0,&\text{then $\wht(v_2)=(0, k)$  and $\wht(f_1)=(0, k-1)$}; \\
2) \ v_1\neq 0,&\text{then $\wht(v_1)=(\ell,1)$ and $\wht(f_1)=(\ell-1,1)$}. \\
\end{array}
\ee

\underline{\textbf{Level 2}: $f_2=\del_1\del_2v\in I(V)$.}
Then, $X^+f_2=0$ implies that $X^+v=0$ and
\[
s_2f_2=(X^-\del_2+\del_1H_1)v=(-\del_1+\del_1H_1+\del_2X^-)v=0
\]
implies that $X^-v=0$ and $H_1v=v$. Since $X^+v=0$, it follows that
\[
\text{$\wht(v)=(1,1)$ and
$\wht(f_2)=(0,0)$. }\qed
\]

\ssec{Proof of Lemma $\ref{L03lev1}$   and Claim $\ref{Cl1}$} In this case, $\fg_0=\fgl(3)$. Denote:
\[
\begin{array}{l}
X_1^+=\xi_1\del_2, \quad X_2^+=\xi_2\del_3, \quad X_3^+=\xi_1\del_3=[X_1^+,X_2^+];\\
X_1^-=\xi_2\del_1, \quad X_2^-=\xi_3\del_2, \quad X_3^-=\xi_3\del_1=-[X_1^-,X_2^-];\\
H_i=\xi_i\del_i\text{~~for~~}i=1,2,3.
\end{array}
\]

We have $\fg_1=\fsvect(0|3)_1\oplus \fsl(1|3)_1$. The subspace $\fsl(1|3)_1$ is spanned by 3 elements:
\[
s_1=\xi_1(\xi_2\del_2+\xi_3\del_3), \; s_2=\xi_2(\xi_1\del_1+\xi_3\del_3), \; s_3=\xi_3(\xi_1\del_1+\xi_2\del_2).
\]
The vector $s_3$ is the lowest weight vector:
\[
[X_2^+,s_3]=s_2, \quad [X_3^+,s_3]=s_1.
\]

The divergence-free component $\fsvect(0|3)_1$ is spanned by 6 elements:
\be\label{u}
\begin{array}{l}
u_1=\xi_2\xi_3\del_1, \; u_2=\xi_3\xi_1\del_2, \; u_3=\xi_1\xi_2\del_3,\\
t_1=\xi_1(\xi_2\del_2-\xi_3\del_3), \; t_2=\xi_2(\xi_3\del_3-\xi_1\del_1), \; t_3=\xi_3(\xi_1\del_1-\xi_2\del_2).
\end{array}
\ee
The vector  $u_1$ is the lowest weight vector:
\[
\ad_{X_1^+}: u_1 \mapsto -t_3\mapsto 2u_2; \; \ad_{ X_3^+}: u_1\mapsto -t_2\mapsto 2u_3; \; \ad_{ X_3^+}: t_3\mapsto -t_1.
\]

Therefore, to check $\fsl(1|3)$-invariance, it suffices to consider only actions of $X_1^+, X_2^+, s_3$, whereas to check $\fvect(0|3)$-invariance, we have to consider, additionally, the $u_1$-action, see~ \eqref{u}.

The brackets we need:
\[
\begin{array}{l}
{}[X_1^+,\del_1]=-\del_2, \; [X_1^+,\del_2]=[X_1^+,\del_3]=0, \; [X_2^+,\del_1]=[X_2^+,\del_3]=0, \; [X_2^+,\del_2]=-\del_3;\\
{}[s_3,\del_1]=-X_3^-,\; [s_3,\del_2]=-X_2^-,\; [s_3,\del_3]=H_1+H_2;\\
{}[u_1,\del_1]=0,\; [u_1,\del_2]=X_3^-, \; [u_1,\del_3]=-X_1^-.
\end{array}
\]

\underline{{\bf Level 1}. Let $f_1=\del_1v_1+\del_2v_2+\del_3v_3\in I(V)$.}

$\bullet$ \textit{Invariance with respect to $\fsl(1|3)$}
\[
\begin{array}{l}
X_1^+f_1=0\Llra \del_1X_1^+v_1-\del_2v_1+\del_2X_1^+v_2+\del_3X_1^+v_3=0.
\\
X_2^+f_1=0\Llra \del_1X_2^+v_1+\del_2X_2^+v_2-\del_3v_2+\del_3X_2^+v_3=0.
\\
s_3f_1=0\Llra -X_3^-v_1-X_2^-v_2+(H_1+H_2)v_3=0.
\end{array}
\]

As a~result, we get the following system of equations:

\vskip 0.3cm

\begin{tabular}{|l|c|c|c|}
\hline
Conditions & $X_1^+f_1=0$ & $X_2^+f_1=0$ & $s_3f_1=0$\\
\hline\hline
Equations & $X_1^+v_1=0$ & $X_2^+v_1=0$ & $X_3^-v_1+X_2^-v_2=(H_1+H_2)v_3$ \\
        & $X_1^+v_2=v_1$ & $X_2^+v_2=0$ &  \\
        & $X_1^+v_3=0$ & $X_2^+v_3=v_2$ &  \\
\hline
\end{tabular}

\vskip 0.3cm

Let the highest weight of $V$ be $(\lambda,\mu, \nu)$.  The following 3 cases are possible:
\be\label{3cases}
\begin{array}{ll}
1)&v_1=v_2=0,\ v_3\ne 0,\\
2)&v_1=0, v_2\ne 0,\\
3)&v_1\ne 0.\\
\end{array}
\ee

$\bullet$ 1) $v_1=v_2=0$, $v_3\ne 0$. Then,
\[
\begin{array}{l}
X_1^+v_3=X_2^+v_3=0 \text{ and } (H_1+H_2)v_3=0\Lra  \\
\wht (v_3)=(n,-n, -n-k) \text{~~and~~}\wht (f_1)=(n,-n, -n-k-1)
\end{array}
\]
where $\begin{cases}n,k\in\Zee_{\geq0}&\text{if $\dim V<\infty$}\\
\text{at least one membership holds: }n,k\in\Cee\setminus \Zee_{\geq0}&\text{if $\dim V=\infty$}.\end{cases}$

If $\dim V=\infty$, the answer is better be written in the form $(\lambda, -\lambda, \nu)$, where at least one membership holds: $\lambda,\nu\in\Cee$.

$\bullet$ 2) $v_1=0$, $v_2\ne 0\Lra X_1^+v_2=X_2^+v_2=0$ and  (taking into account that $X_2^+v_3=v_2$ and $[X_2^+,X_2^-]=H_2-H_3$, $[X_2^+, H_1]=0$, $[X_2^+, H_2]=-X_2^+$)
\begin{equation}\label{v2}
X_2^-v_2=(H_1+H_2)v_3\Lra X_2^+X_2^-v_2=X_2^+(H_1+H_2)v_3\Lra (H_2-H_3)v_2=(H_1+H_2-1)v_2.
\end{equation}

If $\wht(v_2)=(\lambda,\mu,\nu)$, then equation (\ref{v2}) yields
\[
\mu-\nu=\lambda+\mu-1 \Llra \nu=-\lambda+1.
\]

\underline{If $\dim V<\infty$}, then $\lambda-\mu\in\Zee_{\ge 0}$ and $\mu-\nu=\mu+\lambda-1\in\Zee_{\ge 0}$.
We see that
\[
\begin{array}{l}
2\lambda-1\in\Zee_{\ge 0}\Lra \lambda=\nfrac n2,\ \  n\in\Zee_{> 0}; \; \mu=\nfrac n2-k, \ \ k\in\Zee_{\ge 0};\\ \nu=-\nfrac n2+1 \text{ and } \nfrac n2-k+\nfrac n2-1=n-k-1\in\Zee_{\ge 0},
\end{array}
\]
i.e., $0\le k \le n-1$. Finally,
\[
\begin{array}{l}
\wht(v_2)=\left(\nfrac n2, \nfrac n2-k, -\nfrac n2+1\right) \text{~~and~~}\wht (f_1)=\left(\nfrac n2, \nfrac n2-k-1, -\nfrac n2+1\right),\\
\text{ where } 2n\in\Zee_{> 0},\ \  k=0,1,2,\dots, n-1.
\end{array}
\]

\underline{If $\dim V=\infty$}, then $\wht(v_2)=(\lambda, \mu,-\lambda+1)$ and $\wht (f_1)=(\lambda, \mu-1,-\lambda+1)$, where $\lambda-\mu\in\Cee\setminus \Zee_{\ge 0}$.

$\bullet$  3) $v_1\ne 0$. Then,
\[
\begin{array}{l}
X_1^+v_1=X_2^+v_1=0, \; X_2^+v_2=0, \; X_2^+v_3=v_2, \; X_1^+v_2=v_1, \\ (H_1+H_3)v_3=X_3^-v_1+ X_2^-v_2\Lra
X_2^+(H_1+H_3)v_3=X_2^+X_3^-v_1+(H_2-H_3)v_2.
\end{array}
\]
Since
\[
\begin{array}{l}X_2^+(H_1+H_3)v_3=(H_1+H_2-1)X_2^+v_3=(H_1+H_2-1)v_2, \\
X_2^+X_3^-v_1=X_3^-X_2^+v_1+X_1^- v_1=X_1^- v_1,
\end{array}
\]
we get
\[
\begin{array}{l}
(H_1+H_2-1)v_2=X_1^- v_1+(H_2-H_3)v_2\Llra (H_1+H_3-1)v_2=X_1^- v_1\\
\Lra X_1^+(H_1+H_3-1)v_2=X_1^+X_1^- v_1\Lra
\\
(H_1+H_3-2)X_1^+v_2=(H_1-H_2)v_1\Llra (H_1+H_3-2)v_1=(H_1-H_2)v_1\\
\Llra (H_2+H_3-2)v_1=0.
\end{array}
\]

\underline{If $\dim V<\infty$}, then $\mu+\nu=2$ and $\lambda-\mu\in\Zee_{\ge 0}$, \ $2\mu-2\in\Zee_{\ge 0}$. This implies
\[
\mu=\nfrac n2+1 \text{ for }n\in \Zee_{\ge 0}, \ \ \nu=-\nfrac n2+1, \ \ \lambda=\nfrac n2+k \text{ for } k\in\Zee_{> 0}.
\]
Finally,
\be\label{case3}
\begin{array}{l}
\wht(v_1)=\left (\nfrac n2+k, \nfrac n2+1, -\nfrac n2+1\right) \text{~~and~~}\wht (f_1)=\left (\nfrac n2+k-1, \nfrac n2+1, -\nfrac n2+1\right),\\
\text{where } n\in\Zee_{\ge 0},\ \  k\in\Zee_{> 0}.
\end{array}
\ee

\underline{If $\dim V=\infty$}, then $\wht(v_1)=(\lambda, \mu,-\mu+2)$ and $\wht (f_1)=(\lambda-1, \mu,-\mu+2)$, where $\lambda-\mu\in\Cee\setminus \Zee_{\ge 0}$.

$\bullet$ \textit{Invariance with respect to $\fvect(0|3)$}.
We have to add one more condition, $u_1f_1=0$. (Recall, that $u_1v_i=0$.) Since
\[
[u_1,\del_1]=0, \quad   [u_1,\del_2]= X_3^-, \quad  [u_1,\del_3]= -X_1^-,
\]
we get
\begin{equation}\label{v23}
X_3^-v_2=X_1^-v_3.
\end{equation}

$\bullet$
In case 1), i.e., $v_1=v_2=0$, we get an additional condition $X_1^-v_3=0$. Since in this case $X_1^+v_3=0$, we get $(H_1-H_2)v_3=0$. Since in this case we also have $(H_1+H_2)v_3=0$, we conclude that $H_1v_3=H_2v_3=0$, i.e.,
\[
\wht (v_3)=(0,0,-k) \text{~~and $\wht (f_1)=(0,0,-k-1)$, where $k\in\begin{cases}\Zee_{\geq0}&\text{if $\dim V<\infty$},\\
\Cee\setminus \Zee_{\ge 0}&\text{if $\dim V=\infty$}.\end{cases}$}
\]

$\bullet$ In case 2), i.e., $v_1=0, v_2\ne 0$, let us apply $X_3^+$ to (\ref{v23})  and observe that
\[
\begin{array}{l}
X_3^+v_3=(X_1^+X_2^+-X_2^+X_1^+)v_3=X_1^+v_2=v_1;
\\
X_3^+X_3^-v_2=X_3^+X_1^-v_3 \Llra (H_1-H_3)v_2=-X_2^+v_3\\
\Llra (H_1-H_3)v_2=-v_2\Llra (H_1-H_3+1)v_2=0.
\end{array}
\]
Taking into account that $(H_1+H_3-1)v_2=0$ in this case, we conclude: $H_1v_2=0$, $H_3v_2=v_2$.
But in this case $H_1v_2=\lambda v_2$, where $\lambda>0$. Hence,
\[
\begin{cases}\text{no solutions}&\text{if $\dim V<\infty$}\\
\wht(v_2)=(0,\mu,1)\text{~~and $\wht (f_1)=(0,\mu-1,1)$}&\text{if $\dim V=\infty$}.\end{cases}
\]

$\bullet$ In case 3), i.e., $v_1\ne 0$, applying $X_3^+$ to equation (\ref{v23}) we obtain
\[
\begin{array}{l}
(H_1-H_3+1)v_2=X_1^-v_1\Lra X_1^+(H_1-H_3+1)v_2=(H_1-H_2)v_1\Lra
\\
\Lra (H_1-H_3)v_1=(H_1-H_2)v_1\Llra H_3v_1=H_2v_1.
\end{array}
\]
Taking eq. (\ref{case3}) into account, we conclude:
\[
\nfrac n2=-\nfrac n2+2\Lra n=2.
\]
Finally,
\[
\wht(v_1)=(k+1,1,1),  \text{ where $k\in\begin{cases}\Zee_{\geq0}&\text{if $\dim V<\infty$},\\
\Cee\setminus \Zee_{\ge 0}&\text{if $\dim V=\infty$}.\end{cases}$}
\]
 or, in other words,
\[\wht(v_1)=(m,1,1) \text{~~and $\wht (f_1)=(m-1,1, 1)$, where $m\in\begin{cases}\Zee_{>0}&\text{if $\dim V<\infty$},\\
\Cee\setminus \Zee_{> 0}&\text{if $\dim V=\infty$}.\end{cases}$}
\]

\underline{{\bf Level 2}. Let $f_2=\del_1\del_2v_3+\del_1\del_3v_2+\del_2\del_3v_1\in I(V)$.}
In order not to get confused by signs, let us give all the details:
\[
\begin{array}{l}
X_1^+\del_2\del_3=\del_2\del_3X_1^+, \quad X_1^+\del_3\del_1=\del_3\del_1 X_1^+-\del_3\del_2=\del_3\del_1 X_1^++\del_2\del_3, \quad  X_1^+\del_1\del_2=\del_1\del_2 X_1^+,
\\
X_2^+\del_2\del_3=\del_2\del_3X_2^+, \quad X_2^+\del_3\del_1=\del_3\del_1 X_2^+, \quad  X_2^+\del_1\del_2=\del_1\del_2 X_2^+-\del_1\del_3=\del_1\del_2 X_2^++\del_3\del_1,
\\
s_3\del_2\del_3=-\del_2s_3\del_3-X_2^-\del_3=
\del_2\del_3s_3-\del_2(H_1+H_2)-\del_3X_2^-+\del_2
\\
=
\del_2\del_3s_3-\del_2(H_1+H_2-1)-\del_3X_2^-,
\\
s_3\del_3\del_1=-\del_3s_3\del_1+(H_1+H_2)\del_1
=\del_3\del_1s_3+\del_3X_3^-
+\del_1(H_1+H_2-1),
\\
s_3\del_1\del_2=-\del_1s_3\del_2-X_3^-\del_2=\del_1\del_2s_3+\del_1X_2^-
-\del_2X_3^-,
\\
u_1\del_2\del_3=-\del_2u_1\del_3+X_3^-\del_3=
\del_2\del_3u_1+\del_2X_1^-+\del_3X_3^--\del_1,
\\
u_1\del_3\del_1=-\del_3u_1\del_1-X_1^-\del_1=\del_3\del_1u_1-\del_1X_1^-,
\\
u_1\del_1\del_2=-\del_1u_1\del_2=\del_1\del_2u_1-\del_1X_3^-.
\end{array}
\]
Now, consider the actions of $X_1^+, X_2^+, s_3, u_1$ on $f_2$:
\[
\begin{array}{l}
X_1^+f_2=\del_2\del_3X_1^+v_1+\del_3\del_1 X_1^+v_2+\del_2\del_3v_2+\del_1\del_2 X_1^+v_3,
\\
X_2^+f_2=\del_2\del_3X_2^+v_1+\del_3\del_1 X_2^+v_2+\del_1\del_2 X_2^+v_3+\del_3\del_1v_3.
\end{array}
\]
Recall that $s_3v_i=u_1v_i=0$ for $i=1,2,3$. We have
\[
\begin{array}{l}
s_3f_2=-\del_2(H_1+H_2-1)v_1-\del_3X_2^-v_1+\del_3X_3^-v_2+\del_1(H_1+H_2-1)v_2
+\del_1X_2^-v_3-\del_2X_3^-v_3,
\\
u_1f_2=\del_2X_1^-v_1+\del_3X_3^-v_1-\del_1v_1-\del_1X_1^-v_2-\del_1X_3^-v_3.
\end{array}
\]

As a~result, we get the following conditions (the upper line) and the corresponding equations:

\vskip 0.3cm
\be\label{tab}
\begin{tabular}{|c|c|c|c|}
\hline
$X_1^+f_2=0$ & $X_2^+f_2=0$ & $s_3f_2=0$ & $u_1f_2=0$\\
\hline\hline
$X_1^+v_1=-v_2$ & $X_2^+v_1=0$ & $X_2^-v_3=-(H_1+H_2-1)v_2$ &
$X_1^-v_2+X_3^-v_3=-v_1$ \\
$X_1^+v_2=0$ & $X_2^+v_2=-v_3$ & $X_3^-v_3=-(H_1+H_2-1)v_1$ &$X_1^-v_1=0$ \\
$X_1^+v_3=0$ & $X_2^+v_3=0$ & $X_2^-v_1=X_3^-v_2$ &$X_3^-v_1=0$ \\
\hline
\end{tabular}
\ee

\vskip 0.5cm

The following 3 cases are possible.

\vskip 0.3cm

1) $v_2=v_3=0$.

\textit{The $\fsl(1|3)$-invariance}.
In this case, we have:
\[
X_1^+v_1=X_2^+v_2=0, \; (H_1+H_2-1)v_1=0,\; X_2^-v_1=0.
\]
Let us apply $X_2^+$ to the last equation:
\[
X_2^+X_2^-v_1=X_2^-X_2^+v_1+(H_2-H_3)v_1\Lra (H_2-H_3)v_1=0.
\]
We get a~system of equation on the weight of $v_1$:
\[
\wht (v_1)=(\lambda, \mu,\nu)\Lra \left\{
                                  \begin{array}{rcl}
                                  \lambda+\mu-1 & = & 0\\
                                  \mu-\nu & = & 0\\
                                  \end{array}
                                  \right. \Lra \wht (v_1)=(\lambda, -\lambda+1, -\lambda+1).
\]

If $\dim V<\infty$, then $2\lambda-1\in\Zee_{\ge 0}$, i.e.,
\[
\wht(v_1)=\left(\nfrac n2,-\nfrac n2+1, -\nfrac n2+1\right) \text{ for } n\in\Zee_{> 0}.
\]

\textit{The $\fvect(0|3)$-invariance}.
We have to add the condition $u_1f_2=0$. But even the first equation following from this condition yields $v_1=0$. Hence, in this case there are no $\fvect(0|3)$-invariants.

2) $v_3=0, v_2\ne 0$. 
In this case, $X_1^+v_2=X_2^+v_2=0$, and hence $X_3^+v_2=0$, i.e., $v_2$ is the highest weight vector, $\wht(v_2)=(\lambda,\mu,\nu)$.

\textit{The $\fsl(1|3)$-invariance}.
The remaining conditions yield the following system
\begin{equation}\label{level2case2}
    \begin{cases}
    X_1^+v_1 & = -v_2,\\
    X_2^+v_1 & =  0,\\
    (H_1+H_2-1)v_1  & =  0,\\
    (H_1+H_2-1)v_2  & =  0,\\
    X_2^-v_1 & =  X_3^-v_2.\\
    \end{cases}
\end{equation}

First of all, observe that since $[X_1^+,H_1+H_2-1]=0$, the 4th equation of the system (\ref{level2case2}) follows from the 1st and 3rd equations and imposes a~constraint on the weight of  $v_2$:
\[
\lambda+\mu-1=0\Llra \mu=-\lambda +1.
\]

Let us apply $X_3^+$ to the 5th equation and take into account that
\[
[X_3^+, X_2^-]=X_1^+, \text{ and } X_3^+v_1=(X_1^+X_2^+-X_2^+X_1^+)v_1=0.
\]
As a~result, we get:
\begin{equation}\label{nula}
(H_1-H_3)v_2=-v_2 \Llra (H_1-H_3+1)v_2=0\Lra \lambda-\nu+1=0 \Llra \nu=\lambda+1.
\end{equation}

Finally: $\wht(v_2)=(\lambda, -\lambda +1, \lambda +1)$.

If $\dim V<\infty$, then $2\lambda-1\in\Zee_{\ge 0}$ and $-2\lambda\in\Zee_{\ge 0}$, which is impossible.

\textit{ The  $\fvect(0|3)$-invariance}. We have to add the equations corresponding to the condition $u_1f_2=0$:
\begin{equation}\label{l2c2vect}
    \begin{cases}
    X_1^-v_1 & =  0,\\
     X_3^-v_1 & =  0,\\
      X_1^-v_2 & =  -v_1.\\
    \end{cases}
\end{equation}

Applying $X_3^+$ to the 2nd of equations of the system (\ref{l2c2vect}) and taking into account the condition  $X_3^+v_1=0$ obtained earlier, we arrive at the equation $(H_1-H_3)v_1=0$; having applied $X_1^+$ to it we get
\[
X_1^+(H_1-H_3)v_1=(H_1-H_3-1)X_1^+v_1=-(H_1-H_3-1)v_2=0\Lra \lambda-\nu-1=0,
\]
contradicting (\ref{nula}).
\textbf{Verdict}: In this case, there are no $\fvect(0|3)$-invariants. 

\vskip 0.3cm

3) $v_3\ne 0$. 
In this case, $X_1^+v_3=X_2^+v_3=0$,  and hence $X_3^+v_3=0$, i.e., $v_3$ is the highest weight vector, $\wht(v_3)=(\lambda,\mu,\nu)$.

\textit{The $\fsl(1|3)$-invariance}. We have to take into account the condition $s_3f_2=0$ corresponding to the 3 equations from Table \eqref{tab}.

1) We apply $X_2^+$ to eq. $X_2^-v_3=-(H_1+H_2-1)v_2$. The lhs gives
\[
X_2^+X_2^-v_3=X_2^-X_2^+v_3+[X_2^+,X_2^-]v_3=(H_2-H_3)v_3.
\]
The RHS gives:
\[
\begin{array}{l}
X_2^+(-(H_1+H_2-1)v_2)=-((H_1+H_2-1)X_2^+v_2+[X_2^+,H_1+H_2-1]v_2)
\\
=
-(H_1+H_2-2)X_2^+v_2=
-(H_1+H_2-2)(-v_3)=(H_1+H_2-2)v_3.
\end{array}
\]
The second equation is $X_3^-v_3=-(H_1+H_2-1)v_1$ and we apply $X_3^+$. The LHS gives $(H_2-H_3)v_3$. The RHS gives:
\[
\begin{array}{l}
X_3^+(-(H_1+H_2-1)v_1)=-((H_1+H_2-1)X_3^+v_1+[X_3^+,H_1+H_2-1]v_1)\\
=-(H_1+H_2-2)X_3^+v_1.
\end{array}
\]
Let us calculate $X_3^+v_1$:
\[
X_3^+v_1=(X_1^+X_2^+-X_2^+X_1^+)v_1\stackrel{X_2^+v_1=0}{=}-X_2^+(-v_2)=-v_3
\]
which imply a~ condition on the weight of $v_3$:
\[
\lambda+\nu-2=0 \text{ and } \mu+\nu-2=0 \Lra \wht(v_3)=(\lambda,\lambda,-\lambda+2).
\]

It remains to take the 3rd equation into account. Applying $X_3^+$ to it we get
\[
(H_1-H_3)v_2=-X_2^-v_3+X_1^+v_1\Llra(H_1-H_3+1)v_2=-X_2^-v_3.
\]

In more details: We apply $X_3^+$ to $X_3^-v_2=X_2^-v_1$.  Then, the LHS gives
\[
X_3^+X_3^-v_2=X_3^-X_3^+v_2+(H_1-H_3)v_2=(H_1-H_3)v_2
\]
since
\[
X_3^+v_2=(X_1^+X_2^+-X_2^+X_1^+)v_2\stackrel{X_1^+v_2=0}{=}X_1^+(-v_3)=0.
\]
Let us calculate the RHS:
\[
X_3^+X_2^-v_1=X_2^-X_3^+v_1+[X_3^+,X_2^-]v_1\stackrel{X_3^+v_1=-v_3}{=}
-X_2^-v_3+X_1^+v_1.
\]
Having applied  $X_2^+$ to it we get
\[
(H_1-H_3)(-v_3)=-(H_2-H_3)v_3.
\]
As a~result, we arrive at the relation $\lambda=\mu$ obtained earlier.

Finally: $\wht(v_3)=(\lambda,\lambda,-\lambda+2)$.

If $\dim V<\infty$, then $2\lambda-2\in\Zee_{\ge 0}$,  
\[
\wht(v_3)=\left (\nfrac n2+1,\nfrac n2+1,-\nfrac n2+1\right), \text{ where } n\in\Zee_{\ge 0}.
\]

\textit{The $\fvect(0|3)$-invariance}. We have to add the condition $u_1f_2=0$ to the~ 3 equations in Table~ \eqref{tab}.

Having applied $X_3^+$ to the 1st equation, taking into account that $X_3^+v_2=0$, $X_3^+v_1=-v_3$ and $[X_3^+,X_1^-]=-X_2^+$ we get:
\[
X_3^+(X_1^-v_2+X_3^-v_3+v_1)=0\Llra -X_2^+v_2+(H_1-H_3)v_3-v_3=0\Llra (H_1-H_3)v_3=0.
\]

Now, let us apply $X_3^+$ to the 3rd equation twice:
\[
\begin{array}{l}
X_3^+X_3^-v_1=0\Llra X_3^-(-v_3)+(H_1-H_3)v_1=0 \Lra X_3^+X_3^-v_3=X_3^+(H_1-H_3)v_1\Llra
\\
(H_1-H_3)v_3=(H_1-H_3-2)X_3^+v_1\Llra (H_1-H_3-1)v_3=0.
\end{array}
\]

Comparing this result with the previous one, we conclude that $v_3=0$. Thus, in this case there are no $\fvect(0|3)$-invariants.

\textit{Observe that if $\lambda=1$, then the weight of the $\fsl(1|3)$-invariant vector $v_3$ is equal to $(1,1,1)$. However, the vector $f_2$ is not $\fvect(0|3)$-invariant.}

\underline{{\bf Level 3}. Let $f_3=\del_1\del_2\del_3v\in I(V)$.}
Since $\del_i^2=0$ for $i=1,2,3$, we get:
\[
\begin{array}{l}
X_1^+f_3=\del_1\del_2\del_3X_1^+v, \quad X_2^+f_3=\del_1\del_2\del_3X_2^+v,\\
s_3f_3=-\del_1s_3\del_2\del_3v-X_3^-\del_2\del_3v
\\
=-\del_1\del_2\del_3s_3v
+\del_1\del_2(H_1+H_2-1)v+\del_1\del_3X_2^-v-\del_2\del_3X_3^-v+\del_2\del_1v=
\\
=\del_1\del_2(H_1+H_2-2)v-\del_3\del_1X_2^-v-\del_2\del_3X_3^-v;
\\
u_1f_3=-\del_1u_1\del_2\del_3v=-\del_1\del_2X_1^-v-\del_1\del_3X_3^-v.
\end{array}
\]

As a~result, the system of equations describing $\fsl(1|3)$-invariance is as follows.
\begin{equation}\label{level3sl}
\left\{
    \begin{array}{rcl}
    X_1^+v & = & 0, \\
    X_2^+v& = & 0,\\
    (H_1+H_2-2)v& = & 0,\\
    X_2^-v& = & 0, \\
    X_3^-v& = & 0. \\
    \end{array}
\right.
\end{equation}
The first two of equations \eqref{level3sl} imply $X_3^+v=0$. Applying $X_2^+$ to the 4th equation, and $X_3^+$ to the 5th, we get
\begin{equation}\label{lev3wht}
(H_2-H_3)v=0 \text{ and } (H_1-H_3)v=0,
\end{equation}
implying that $\wht(v)=(\lambda,\mu,\nu)$ satisfies the condition $\lambda=\mu=\nu$. The 3rd of equations \eqref{level3sl} implies $\lambda+\mu=2$, so finally we have
\[
\wht(v)=(1,1,1). \hfill \qed
\]

\ssec{Proof of Lemma $\ref{L20}$} For a~basis of $\fg_0$, where $\fg=\fvect(2|0)$ in the standard grading, we take $X^-:=x_2\del_1$, $H_i:=x_i\del_i$ and $X^+:=x_1\del_2$ and for a~basis of $\fg_1$ we take $s_i:=x_iE$, where $E:=\sum x_i\del_i$, and the 4 divergence-free elements
spanning an irreducible $\fgl(2)$-module which are not needed to determine Bols. We will need the following relations
\be\label{ut}
\begin{array}{l}
{}[X^+,\del_1]=-\del_2, \  [[X^+,\del_1],\del_1]=0 \Lra X^+\del^n_1=\del^n_1X^+-n\del^{n-1}_1\del_2,\\
{}[X^+,\del_2]=0\Lra X^+\del^n_2=\del^n_2X^+ \Lra X^+\del^m_1\del^k_2=\del^m_1\del^k_2X^+-m\del^{m-1}_1\del_2^{k+1},\\
{}[s_2,\del_1]=-X^-, \; [[s_2,\del_1],\del_1]=0 \Lra s_2\del_1^n=\del_1^ns_2-n
\del_1^{n-1}X^-,\\
{}[s_2,\del_2]=-H_1-2H_2.
\end{array}
\ee

Auxiliary identities:
\begin{equation}\label{H12}
\begin{array}{l}
{}[H_1+2H_2,\del_2]=-2\del_2,\; [[H_1+2H_2,\del_2],\del_2]=0\\
 \Lra (H_1+2H_2)\del_2^n=\del_2^n(H_1+2H_2)-2n\del_2^n.
 \end{array}
\end{equation}

\begin{equation}\label{x-}
[X^-,\del_2]=-\del_1, \; [[X^-,\del_2],\del_2]=0 \Lra X^-\del_2^k=\del_2^kX^- -k\del_1\del_2^{k-1}.
\end{equation}

Equation (\ref{H12}) implies that
\[
s_2\del_2^n=\del_2^ns_2-n\del_2^{n-1}(H_1+2H_2)+a_n\del_2^{n-1}.
\]
Hence,
\[
\begin{array}{l}
s_2\del_2^{n+1}=\del_2s_2\del_2^n-(H_1+2H_2)\del_2^n\\
=\del_2^{n+1}s_2-n\del_2^n(H_1+2H_2) +a_n\del_2^n-\del_2^n(H_1+2H_2)+2n\del_2^2\\
=\del_2^{n+1}s_2-(n+1)\del_2^n(H_1+2H_2)+(a_n+2n)\del_2^n\Lra a_{n+1}=a_n+2n.
\end{array}
\]

Since $a_1=0$ and $a_2=2$, we deduce that $a_n=n(n-1)$. Therefore,
\[
\begin{array}{ll} 
s_2\del_2^n&=\del_2^ns_2-n\del_2^{n-1}(H_1+2H_2)+n(n-1)\del_2^{n-1}\\
&=
\del_2^ns_2-n\del_2^{n-1}(H_1+2H_2-n+1).
\end{array}
\]

Finally, with eq. (\ref{x-}) we get
\[
\begin{array}{l}
s_2\del_1^m\del_2^k=\del_1^ms_2\del_2^k-m\del_1^{m-1}X^-\del_2^k\\
=
\del_1^m\del_2^ks_2-k\del_1^m\del_2^{k-1}(H_1+H_2-k+1)-m\del_1^{m-1}\del_2^kX^-
+mk\del_1^m\del_2^{k-1}\\
=\del_1^m\del_2^ks_2-k\del_1^m\del_2^{k-1}(H_1+H_2-m-k+1)-m\del_1^{m-1}\del_2^kX^-.
\end{array}
\]

Summary of the above:
\[
\begin{array}{lcl}
X^+\del^n_1  & = & \del^n_1X^+-n\del^{n-1}_1\del_2,\\
X^+\del^n_2  & = & \del^n_2X^+, \\
X^+\del_1^m\del^k_2 & =  & \del^m_1\del^k_2X^+-m\del^{m-1}_1\del_2^{k+1},\\
s_2\del_1^n & = & \del_1^ns_2-n\del_1^{n-1}X^-,\\
s_2\del_2^n & = & \del_2^ns_2-n\del_2^{n-1}(H_1+2H_2-n+1),\\
s_2\del_1^m\del_2^k & = & \del_1^m\del_2^ks_2-k\del_1^m\del_2^{k-1}(H_1+H_2-m-k+1)-m\del_1^{m-1}\del_2^kX^-.
\end{array}
\]

{\bf Level $n$:} Let $f_n=\sum_{m=0}^n\del_1^m\del_2^{n-m}v_m$ be a singular vector. Then,
\[
X^+f_n=\del_1^nX^+v_n + \sum_{0\leq m\leq n-1} \del_1^m\del_2^{n-m}(X^+v_m-(m+1)v_{m+1}),
\]
implying
\be\label{SystX+}
X^+f_n=0 \Llra \left\{\begin{array}{lcl}
                      X^+v_n & = & 0,\\
                      X^+v_m & = & (m+1)v_{m+1}\text{~~for~~} m=  0,1,\dots, n-1.\\
                      \end{array}
                \right.  
\ee
Since $s_2v_m=0$ for $m=0,1,\dots, n$,  the condition $s_2f_n=0$ is equivalent to the system
\be\label{Syst}
 mX^-v_m+(n-m+1)(H_1+2H_2-n+1)v_{m-1}=0 \text{ for all } m=1,2,\dots n.
\ee
Let us solve eq.~\eqref{Syst}. 
There are $n+1$ cases:
\[
\textbf{Case $k$: } v_n=\dots = v_{k+1}=0, \ \ v_k\ne 0, \text{ where } k= 0,1,\dots, n.
\]

Let  $\wht(v_k)=(\lambda, \mu)$ denote the highest weight of $V$. Set $H:=H_1+2H_2-n+1$. 

{\bf Case 0}. 
Then, conditions (\ref{SystX+}) and (\ref{Syst}) yield just one condition:
\[
Hv_0=0 \Lra \wht (v_0)=(n-1-2\lambda, \lambda), \ \ \wht (f_n)=(n-1-2\lambda, \lambda-n).
\]

{\bf Case $k$, where $0< k<n$ and $n\ge 2$}. 
In this case, each of the systems  (\ref{SystX+}) and (\ref{Syst}) contains no fewer than 3 equations. In more details, where $m<k$:
\[
\left\{  \begin{array}{rclcc}
           X^+v_k & = &  0 & & (X_k)\\
           X^+v_{k-1} & = & k v_k & & (X_{k-1})\\
           \dots & \dots & \dots &&\\
           X^+v_{m} & = & (m+1) v_{m+1} & & (X_m)\\
           \dots & \dots & \dots &&\\
           X^+v_0 & = & v_1 & & (X_0)\\
           Hv_k & = & 0 && (S_{k+1})\\
           kX^-v_k +(n-k+1)Hv_{k-1} & = & 0 && (S_k)\\
           \dots & \dots & \dots &&\\
           (m+1)X^-v_{m+1} +(n-m)Hv_{m} & = & 0 && (S_{m+1})\\
            mX^-v_{m} +(n-m+1)Hv_{m-1} & = & 0 && (S_m)\\
           \dots & \dots & \dots &&\\ 
            X^-v_1+nHv_0 & = & 0 && (S_1)\\
          \end{array}
\right.
\]

We have to understand what constraints on $\wht(v_k)$ do equations $(S_1)-(S_{k+1})$ impose. 

Applying $X^+$  to $(S_m)$ we get: 
\[
m(m+1)X^-v_{m+1}+m(H_1-H_2)v_m+m(n-m+1)(H+1)v_m=0.
\]

Brake $m$ out and subtract $(S_{m+1})$:
\be\label{Bez}
(H_1-H_2+H+(n-m+1))v_m=0.
\ee

Apply $(X^+)^{k-m}$ to eq. \eqref{Bez}, take into account that $[X^+,H_1-H_2+H]=-X^+$ and use equations $(X_m)-(X_{k-1})$; we get
\begin{equation}\label{Sm}
(H_1-H_2+H+(n-k+1))\frac{k!}{m!}v_k=0 \Llra (H_1-H_2+H+(n-k+1))v_k=0.
\end{equation}

This equation does not depend on $m$, i.e., all equations $(S_1)-(S_{k-1})$ yields just one equation~ (\ref{Sm}). As a result, we are left with just 3 equations:
\[
\left\{  \begin{array}{rclcc}
           Hv_k & = & 0 && (S_{k+1})\\
           kX^-v_k +(n-k+1)Hv_{k-1} & = & 0 && (S_k)\\
           (H_1-H_2+H+(n-k+1))v_k& = &0 && (S_{<k})\\
          \end{array}
\right.
\]
Now, apply $X^+$ to equation $(S_k)$ and take $(X_{k})$ into account; we get
\[
\begin{array}{l}
 k(H_1-H_2)v_k+k(n-k+1)(H+1)v_k=0  \\
 \Llra (H_1-H_2+H+n-k+1)v_k+(n-k)Hv_k =0, 
 \end{array}
\]
i.e., equation $(S_{k+1})$ follows from $(S_{k})$ and $(S_{<k})$.

As a result, we get a system of equations on the weight $\wht(v_k)=(\lambda,\mu)$:
\[
\left\{  \begin{array}{rcl}
          \lambda+2\mu & = & n-1,\\
          2\lambda+\mu & = & k-2\\  
          \end{array}
\right.  \Llra 
\left\{  \begin{array}{rcl}
          \lambda & = & \frac{2k-n}3-1,\\
          \mu & = & \frac{2n-k}3.\\
          \end{array}
\right.
\]
Then, $\dim V=\infty$ and
\[
\wht(f_n)=(-\nfrac{k+n}3-1, \nfrac{2k-n}3).
\]

{\bf Case $n$}. 
Unlike the previous case, there is no condition
$(S_{n+1})$, 
we only have 
\[
\left\{  \begin{array}{rclcc}
           nX^-v_n +Hv_{n-1} & = & 0 && (S_n)\\
           (H_1-H_2+H+1)v_n& = &0 && (S_{<n})\\
          \end{array}
\right.
\]
Applying $X^+$ to equation $(S_n)$ and taking  $(X_{n})$ into account we get
\[
 n(H_1-H_2)v_k+n(H+1)v_n=0  \Llra (H_1-H_2+H+1)v_n =0,
\]
i.e., precisely $(S_{<n})$. Thus, in this case, we have only one equation on $\wht(v_n)=(\lambda,\mu)$:  
\[
2\lambda+\mu=n-2 ,
\]
and therefore 
\[
\wht(v_n)=(\lambda, n-2-2\lambda), \ \ \wht(f_n)=(\lambda-n, n-2-2\lambda). 
\]
Note that for $n=1$ we have only one equation $(S_1)$ from the very beginning. \qed

\ssec{Proof of Lemma $\ref{L11}$ and Claim $\ref{C11}$}  For a~basis of $\fg_{0}$ we take $X^-:=\xi\del$, $H_1:=x\del$, $H_2:=\xi\delta$ and $X^+:=x\delta$;  for a~basis of $\fg_{1}$ we take $s_x:=xE$ and $s_\xi:=\xi E=x\xi\del$, where $E=x\del+\xi\delta$, and divergence-free elements $u_1:=x^2\del+2x\xi\delta$ and $u_2:=x^2\delta$ of which the lowest weight vectors are $s_\xi$ and  $u_1$. Then,
\[
\begin{array}{l}
X^+(f_n)=-n\del^{n-1}\delta v+\del^nX^+v-\del^{n-1}\delta X^+ w =0\Lra X^+w=-nv,\ \ X^+v=0;\\
s_\xi(f_n)=-\del^{n-1}(nX^-v+(n-1-H_1)w) \Lra nX^-v+(n-1-H_1)w=0.
\end{array}
\]
Take into account that 
\[
\begin{array}{l}
X^+(nX^-v+(n-1-H_1)w)=n(H_1+H_2)v+(n-1)X^+w+X^+w-H_1X^+w\\
\Lra (H_1+H_2)v-(n-1)v-v+H_1v=(2H_1+H_2-n)v.
\end{array}
\]

Now, let us describe the singular vectors corresponding to the $\fvect(1|1)$-invariant differential operators. Consider one more condition
\be\label{v}
u_1(f_1)= -2(H_1+H_2)v-2X^+w=-2(H_1+H_2-1)v.
\ee
There are two cases, where $\lambda\in\Cee$:
\[
\begin{array}{ll}
1) \ v=0,&\text{then~ no new conditions, same answer as in Lemma $\ref{L11}$};\\
2) \ v\neq 0,&\text{then~}  \wht(v)=(\lambda, 1-2\lambda),\ \ 1-2\lambda=-\lambda+1 \Lra \lambda=0, \wht(f_1)=(-1,1).  \\
\end{array}
\]

Direct computation shows that the condition $u_1(f_2)=0$ leads to a~contradiction, therefore there are no singular vectors of level $>1$. \qed

 \subsection*{Acknowledgments}
S.B. and D.L. were supported by the grant NYUAD 065.

\def\eightit{\it}
\def\bib{\bf}
\bibliographystyle{amsalpha}

\end{document}